\documentclass[11pt]{amsart}
\usepackage{amsmath,amssymb,graphicx}
\usepackage{amssymb}
\usepackage{pb-diagram}
\usepackage{verbatim}

\newtheorem{thm}{Theorem}
\newtheorem{lem}[thm]{Lemma}

\newtheorem{prop}[thm]{Proposition}

\theoremstyle{definition}
\newtheorem{defn}[thm]{Definition}

\newtheorem{rmk}[thm]{Remark}
\newtheorem{exmp}[thm]{Example}

\newcommand{\CPb}{\overline{\mathbb{CP}}{}^{2}}
\newcommand{\CP}{{\mathbb{CP}}{}^{2}}
\newcommand{\R}{\mathbb{R}}
\newcommand{\Z}{\mathbb{Z}}

\title[Simply Connected Symplectic Calabi-Yau $6$-Manifolds]
{Simply Connected Symplectic Calabi-Yau $6$-Manifolds}

\begin{document}

\author{Anar Akhmedov}
\address{School of Mathematics,
University of Minnesota,
Minneapolis, MN, 55455, USA}
\email{akhmedov@math.umn.edu}

\date{June 03, 2011.  Revised on July 10, 2011}

\subjclass[2000]{Primary 57R55; Secondary 57R17}

\begin{abstract} In this article, we construct simply connected symplectic Calabi-Yau $6$-manifold by applying Gompf's symplectic fiber sum operation along $\mathbb{T}^4$. Using our method, we also construct symplectic non-K\"{a}hler Calabi-Yau $6$-manifolds with fundamental group $\Z$. We also produce the first examples of simply connected symplectic Calabi-Yau and non-Calabi-Yau $6$-manifolds via coisotropic Luttinger surgery on non simply connected symplectic $6$-manifolds.
\end{abstract}

\maketitle
\section{Introduction}

Our work in this paper is inspired and motivated by works of I. Smith, R. Thomas, and S. T. Yau \cite{STY}, and more recent works of T-J. Li and C-I. Ho \cite{Ho, HoLi}, and S. Baldridge and P. Kirk \cite{BK} where they introduced the symplectic surgery operations on symplectic $6$-manifolds, called the symplectic conifold transition and coisotropic Luttinger surgery, respectively to study the symplectic Calabi-Yau $6$-manifolds. The main goal of this paper is to construct simply connected symplectic Calabi-Yau $6$-manifolds. Our construction will employ Gompf's symplectic connected sum operation. Along the way, we also produce symplectic non-K\"{a}hler Calabi-Yau $6$-manifolds with fundamental group $\Z$. The only known purely symplectic construction of simply connected symplectic Calabi-Yau $6$-manifold were given recently by J. Fine and D. Panov in their seminal paper \cite{FP2} in 2009. We believe that our construction is simpler than than the construction in \cite{FP2}. The examples constructed more recently by Baldridge-Kirk in \cite{BK} via coisotropic Luttinger surgery have $b_1 \geq 2$. Very little is known about the geography of symplectic Calabi-Yau $6$-manifolds. A concise history of constructing symplectic Calabi-Yau $6$-manifolds can be found in \cite{STY, FP1, FP2, BK, Ho}. In comparision, we would like to remark that, in dimension $4$, a great deal is known about the geography of symplectic Calabi-Yau manifolds (\cite{MS}, \cite{Li1}, \cite{Li3}, \cite{B})(see also the survey \cite{Li2}). 

Let $\CP$ denote the complex projective plane and let $\CPb$ denote the underlying smooth $4$-manifold $\CP$ equipped with the opposite orientation. Let $X$ denote $E(1) = \CP\# 9\CPb$, the complex projective plane blown up at $9$ points. Our main results are the following theorems.

\begin{thm}\label{thm:main1} There exist a simply connected symplectic\/ Calabi-Yau $6$-manifold that can be obtained from $E(1) \times \mathbb{T}^2$ by symplectic connected sum along $\mathbb{T}^4$.
\end{thm}

\begin{thm}\label{thm:main2} There exist a non-K\"{a}hler Calabi-Yau symplectic\/ $6$-manifolds with the fundamental group $\Z$ that can be obtained from $E(1) \times \mathbb{T}^2$ by symplectic connected sum along $\mathbb{T}^4$. \end{thm}

Our paper is organized as follows. Section~\ref{sec: SCS} contains a brief review of symplectic surgery operation on symplectic manifolds called the symplectic connected sum, and the proofs of some preliminary results that will be used in our proof of Theorem~\ref{thm:main1} and Theorem~\ref{thm:main2}. In Section~\ref{sec: SPCY}, we construct simply connected symplectic Calabi-Yau $6$-manifolds as the symplectic connected sum of two copies of symplectic $6$-manifold $E(1) \times \mathbb{T}^2$ along $\mathbb{T}^4$ and present proof of our main Theorem~\ref{thm:main1}. In Section~\ref{sec: FCY}, we construct symplectic non-K\"{a}hler Calabi-Yau $6$-manifolds with fundamental group $\Z$ and prove Theorem~\ref{thm:main2}. Finally, in Section~\ref{CLS}, we also show how to obtain such examples (simply connected and $\pi_1 = \Z$ examples) by coisotropic Luttinger surgery on many non-simply connected symplectic $6$-manifolds. Several other applications are also discussed in Section~\ref{CLS}. In a forthcoming paper, we study in detail the geography of symplectic Calabi-Yau $6$-manifolds \cite{AL}.

\section{Symplectic Connected Sum}{\label{sec: SCS}} In this Section, we recall some basic facts about the symplectic connected operation, see \cite{gompf} for more details. In \cite{gompf} Gompf defines a surgery operation which takes two smooth manifolds and glues them along diffeomorphic submanifolds. Do do such gluing, one needs an orientation reversing diffeomorphism of the normal bundles of the given submanifolds. Furthermore, he has proved that this is a symplectic operation. Namely, if the manifolds and submanifolds are symplectic, then the resulting manifold will admit a symplectic structure as well.

\begin{defn} Let $(X_1, \ \omega_{X_1})$ and $(X_2, \ \omega_{X_2})$ be closed symplectic $2n$-dimensional manifolds. Suppose that $(Y, \ \omega_{Y})$ is another symplectic manifold of dimension $2n-2$, and that there exists symplectic embeddings $j_{i}: Y \rightarrow X_i$. Denote by $Y_i$ the images $j_i(Y)$ and by $\nu_i$ their normal bundles in $X_i$. Assume that the Euler class of $\nu_i$ satisfy $e(\nu_1) +  e(\nu_2) = 0$. Then for any choice of an orientation reversing $\psi: \nu_1 \cong \nu_2$, the \emph{symplectic connected sum} of $X_1$ and $X_2$ along $Y$ is the manifold $(X_1\setminus Y_1)\cup_{\psi} (X_2\setminus Y_2)$ and is denoted by $X_1 \#_{\psi} X_2$

\end{defn}

The diffeomorphism type of $X_1 \#_{\psi} X_2$ depends on the choice of the embeddings and of the orientation reversing bundle isomorphism $\psi$. We will sometimes abuse our notation and denote the symplectic connected sum by $X_1 \#_{Y} X_2$. 

\begin{thm}\label{thm:symsum} For any choice of an orientation reversing $\psi: \nu_1 \cong \nu_2$, the manifold $X_1 \#_{\psi} X_2$ admits a canonical symplectic structure $\omega$ induced by $\omega_{X_1}$ and $\omega_{X_2}$.

\end{thm}

Consider the union of $X_i$ along $Y_i$. One can construct a fibration $Z$ over $D^2$ with this union as the central fiber and the sum manifold as a general fiber. In fact in \cite{IP} it is shown that there is a symplectic form $\Omega$ on $Z$ such that it restricts to $\omega$ on the general fiber, and to $\omega_i$ on $X_i$. For the convenience of the reader, we state E. Ionel and T. Parker's result below. 

\begin{thm}\label{thm:IP} There exists a $2n + 2$-dimensional symplectic manifold $(Z, \Omega)$ and a fibration $\lambda: Z \rightarrow D$ over a disk $D \subset \mathbb{C}$. The central fiber $Z_0$ is the singular symplectic manifold $X_1 {\cup}_Y X_2$ , while for $\lambda \ne 0$,
the fibers $Z_{\lambda}$ are smooth compact symplectic submanifolds - the symplectic
connect sums.
\end{thm}

In \cite{IP} E. Ionel and T. Parker also prove a useful lemma comparing the first Chern class of the symplectic sum $X = X_1 \#_{\psi} X_2$ (for arbitrary even dimension) with the Chern classes $c_1(X_1)$ and $c_1(X_2)$ of $X_1$ and $X_2$. 

\begin{lem}\label{lem:IP} If $A \in H_2(Z_{\lambda}; \Z)$, $\lambda \ne 0$, is homologus in $Z$ to the union  $C_1\cup C_2$ in $X_1\cup_Y X_2$ of cycles $C_i$ in $X_i$, then
\begin{equation}
\begin{array}{l}
c_1(Z)(A) = c_1(Z_{\lambda})(A) = c_1(X)(A) = \\
c_1(X_1)([C_1])+c_1(X_2)([C_2])-\sum [Y_i]\cdot [C_i].
\end{array}
\end{equation}
and for a rim class $r$, $c_1(X)(r)=0$.
\end{lem}

We will need the following proposition in our proofs of Theorems \ref{thm:main1} and \ref{thm:main2}.

\begin{prop}\label{prop:symsum} Let $X$ be closed, symplectic manifold of dimension $2n$, and $Y$ be closed, symplectic submanifold of of dimension $2n-2$ with trivial normal bundle. If there exist a sphere $S$ in $X$ that intersects $Y$ transversally in exactly one point, then the homomorphism $j_{*}: \pi_1(X \setminus Y) \rightarrow \pi_1(X)$ induced by inclusion is an isomorphism. In particular, if $X$ is simply connected, then so is $X \setminus Y$.

\end{prop} 

\begin{proof}

We closely follow the proof of Proposition 1.2 in \cite{Halic}, see also Gompf's article \cite{gompf} where this criterion was initially used. Notice that $Y$ has a symplectically embedded neighborhood $Y \times D_{\epsilon}$ in $X$. Since $Y$ has codimension $2$ in $X$, any loop in $X$ can be homotoped away from $Y$, thus the homomorphism $j_{*}: \pi_1(X \setminus Y) \rightarrow \pi_1(X)$ is surjective. Let $\gamma$ be any loop in $\pi_1(X \setminus Y)$ such that the image $j_{*}[\gamma] = 0$. We consider a homotopy $H: I \times I \rightarrow X$ of $\gamma$ to a constant map. The homotopy $H$ can be choosen in a way that it meets $Y \times \partial(D_{\epsilon})$ in a finite number of circles $\gamma_{k}$, thus we have $[\gamma] = \prod_{k} [\gamma_k]$ in $\pi_1(X \setminus Y)$. To show the homotopy class $[\gamma]$ is trivial in $\pi_1(X \setminus Y)$, it suffices to show the homotopy classes $[\gamma_k]$ are all trivial in $\pi_1(X \setminus Y)$. We can move each circle $\gamma_k$ until we reach the intersection circle $\alpha : = S \cap (Y \times \partial(D_{\epsilon}))$. Since $\alpha$ is null-homotopic in $S \setminus {point} = \mathbb{C}$, so are the circles $\gamma_{k}$. This shows that $j_{*}: \pi_1(X \setminus Y) \rightarrow \pi_1(X)$ induces an isomorphism.          

\end{proof}

\begin{defn} A \emph{symplectic Calabi-Yau manifold}, CY for short, is a symplectic manifold $M$ with $c_1(M) = 0$.
\end{defn}

\begin{exmp}\label{E} Let $E(n)$ be a simply connected elliptic surface without multiple fibers. In this example, we study the elliptic surfaces $E(2) =  E(1) \#_{\mathbb{T}^2} E(1)$, $K3$ surface, in some detail. Our discussion will be useful in our proof of the Theorems \ref{thm:main1} and \ref{thm:main2}. We will think of $K3$ surface as the fiber sum of two copies of $E(1) = \CP\# 9\CPb$ along a torus fiber. Consider the following basis for the intesection form of $E(1)$:  $ <f = 3h - e_1 - ... - e_9, \ e_9, \ e_1 - e_2, \ e_2 - e_3, \ \cdots, \ e_7 - e_8, \ -h + e_6 + e_7 + e_8 >$, where $e_i$ denote the homology class of the exceptional sphere of the $i-th$ blow up and $h$ the pullback of the hyperplane class of $\CP$. The last $8$ classes can be represented by spheres of self-intersection $-2$ and generate the intersection matrix  $-E_8$, where $E_8$ the matrix corresponding to the Dynkin diagram of the exceptional Lie algebra $E_8$. The class $f$ is fiber of an elliptic fibration on $E(1) = \CP\# 9\CPb$ and $e_9$ is a section. When we perform the fiber sum to get $E(2)$, it is not hard to see the surfaces that generate the intersection form $2(-E_8) \oplus 3H$ for $E(2)$, where $H$ is a hyperbolic pair. The two copies of the Milnor fiber $\Phi(1) \in E(1)$ are in $E(2)$, providing $16$ spheres of self-intersection $-2$ (corresponding to the classes $\{ e_1 - e_2, \ e_2 - e_3, \ \cdots, \ e_7 - e_8, \ -h + e_6 + e_7 + e_8 \}$ mentioned above), which realize two copies of $-E_8$. One copy of $H$ comes from a torus fiber $f$ and a sphere section $\sigma$ of self-intersection $-2$, i.e. from the Gompf's nucleus $N(2)$ in $E(2)$. The remaining two copies of $H$ come from $2$ rim tori and their dual $-2$ spheres (see discussion in \cite{GS}, page 73)). These $22$ classes ($19$ spheres and $3$ tori) generate $H_2$ of $E(2)$. Since $c_1(E(n)) = (2-n)f$, $E(2)$ is CY manifold.  

\end{exmp}

\begin{prop}  The symplectic connected sum $X_1\#_{\psi}X_2$ of $X_1$ and $X_2$ along $Y$ has Chern numbers given by $c_I [X_1\#_{\psi}X_2] = c_I [X_1]+c_I[X_2] - c_I [Y \times S^2]$ where $I$ stands for any arbitrary partition of $n$.
\end{prop}

For the proof of above Proposition, we refer the reader to \cite{Halic}. 

In the discussion that follows, we will only consider the symplectic $6$-manifolds. From the above proposition, we easily obtain the following formulas for the Chern numbers of the symplectic connected sum $X_1\#_{\psi} X_2$:

\begin{equation}\label{eq: Chern Numbers}
\begin{array}{l}
{c_1}^3(X_1\#_{\psi}X_2) = {c_1}^3(X_1) + {c_1}^3(X_2) - 6{c_1}^2(Y), \\
{c_1c_2}(X_1\#_{\psi}X_2) = {c_1c_2}(X_1) + {c_1c_2}(X_2)- 2({(c_1}^2(Y) + c_2(Y)), \\
{c_3}(X_1\#_{\psi}X_2) = c_3(X_1) + c_3(X_2) - 2c_2(Y)
\end{array}
\end{equation}

\section{Proof of Theorem~\ref{thm:main1}}
{\label{sec: SPCY}}


\begin{proof}

To construct our simply connected symplectic Calabi-Yau $6$-manifolds, we take two copies $W = E(1) \times T^2$ and form their symplectic fiber sum along the symplectic submanifolds $F \times \mathbb{T}^2$ and $F' \times \mathbb{T}^2$, where $F$ and $F'$ are the regular fibers of an elliptic fibration on $E(1)$. Notice that $-1$ sphere section $S$ of $E(1)$ gives a section for the fibration on $E(1) \times \mathbb{T}^2$ with a regular fiber $F \times \mathbb{T}^2 = \mathbb{T}^4$. By Proposition \ref{prop:symsum}, $\pi_1(W \setminus F \times \mathbb{T}^2) = \pi_1(W) = \Z \times \Z$. The fundamental group of $W$ is generated by the circles $c$ and $d$ coming from product $\mathbb{T}^2$, the generators $a$ and $b$ of $F$ and the normal circle $\mu$ are all nullhomotopic in $\pi_1(W \setminus F \times \mathbb{T}^2)$. Our symplectic $6$-manifold will be the symplectic connected sum of two copies $W$ along the $4$-tori $F \times \mathbb{T}^2$ and $F' \times \mathbb{T}^2$, i.e. $X_{\psi} = W \#_{\psi} W$. Let us choose a special gluing diffeomorphism $\psi : \partial(F \times \mathbb{T}^2 \times D^2) \longrightarrow  \partial(F' \times \mathbb{T}^2 \times D^2)$ that comes from an orientation preserving diffemorphism  of $\mathbb{T}^4$ which sends the generators of $\pi_1$ as follows:
\begin{equation*}
a = 1 \mapsto c', \ \
b = 1 \mapsto d', \ \
c \mapsto a' = 1, \ \
d \mapsto b' = 1 , \ \
\mu = 1 \mapsto {\mu'}^{-1} = 1, \\
\end{equation*}

\vspace{2pt}

From Seifert-Van Kampen Theorem, we get the following presentation for the fundamental group of $X_{\psi}$.

\begin{eqnarray}
\pi_1(X_{\psi})  &=& \langle
c, d;\ c', d' \mid c=1,\, d=1,\, c'=1,\, d'=1\rangle
\end{eqnarray}

\vspace{2pt}

This shows that $\pi_1(X_{\psi}) = 1$. By Gompf's Theorem \ref{thm:symsum}, $X_{\psi}$ is symplectic.

\vspace{3pt}

Let us now prove that the simply connected symplectic $6$-manifold $X_{\psi}$ is Calabi-Yau. We need to verify that $c_1(X_{\psi}) = 0$. 

\vspace{2pt}

First, we construct a basis for $2$-dimensional homology of $X_{\psi}$. Since $X_{\psi}$ is obtained as the symplectic connected sum of two copies of $E(1) \times T^2$, it is relatively easy to construct such a basis using Example ~\ref{E}. To get a basis for $H_2$ of $E(1) \times T^2$ ($b_2 = 11$), we use the following basis for the two dimensional homology of $E(1)$:  $ <f = 3h - e_1 - ... - e_9, \ e_9, \ e_1 - e_2, \ e_2 - e_3, \ \cdots, \ e_7 - e_8, \ -h + e_6 + e_7 + e_8 >$ in Example ~\ref{E} and $pt \times T^2$. When we do the symplectic connected sum along $F \times \mathbb{T}^2$ to get our symplectic $6$-manifold $X_\psi$, because of our choice of $\psi$, the rim tori and their associated vanishing classes are all nullhomologous in $X_\psi$. The reason for that is the following: the circles $c$, $d$, $c'$, and $d'$ in $\partial(W \setminus \mathbb{T}^2 \times \mathbb{T}^2 \times D^2)$ do not bound disks on both sides. Thus, a basis for $2$-dimensional homology of $X_{\psi}$ comes from the two copies of the Milnor fiber $\Phi(1) \in E(1) \subset E(1) \times T^2$, $-2$ sphere section $\sigma'$ obtained by sewing the $-1$ sphere sections $e_9$ and ${e_9}'$, the fiber $f$, and $pt \times T^2$. Let us show that $c_1(X_\psi)$ is zero on all of these two dimensional classes. We will use the lemma~\ref{lem:IP} of Ionel and Parker above. Since $c_1(E(1)\times \mathbb{T}^2) = PD(F \times \mathbb{T}^2)$, $c_1(X_\psi)$ evaluates to zero on $2$-dimensional homology classes of $X_{\psi}$ coming from the two copies of the Milnor fiber $\Phi(1)$. Notice that these $16$ spheres have no intersection with $F \times \mathbb{T}^2$. $c_1(X_\psi)$ is zero on sphere section $\sigma'$ coming from the union of $-1$ sphere sections of $E(1)$, thus zero by Lemma \ref{lem:IP}. Since the push-offs of the surfaces $f$ and $pt \times \mathbb{T}^2$ in the normal direction of $F \times \mathbb{T}^2$ have no intersection with $F \times \mathbb{T}^2$, $c_1(X_\psi)$ is zero on these classes. Since $c_1(X_\psi)$ is zero on our basis, we have $c_1(X_\psi)$ is zero class.

\vspace{2pt}

This concludes the proof of the Theorem.

\end{proof}

Using the formulas above, setting $X_1 = X_2 = E(1) \times \mathbb{T}^2$, $Y = \mathbb{T}^4$, and ${c_1}^2(\mathbb{T}^4) = c_2(\mathbb{T}^4) = {c_1}^3( E(1) \times \mathbb{T}^2) = {c_1c_2}(E(1) \times \mathbb{T}^2) = 0$,
we can easily compute the Chern numers of $X_{\psi}$.

\begin{equation}\label{eq: Chern Numbers X}
\begin{array}{l}
{c_1}^3(X_{\psi}) = 2{c_1}^3(E(1)\times \mathbb{T}^2) - 6{c_1}^2(\mathbb{T}^4) = 0, \\
{c_1c_2}(X_{\psi}) = 2{c_1c_2}(E(1)\times \mathbb{T}^2) - 2({(c_1}^2(\mathbb{T}^4) + c_2(\mathbb{T}^4)) = 0, \\
{c_3}(X_{\psi}) = 2c_3(E(1)\times \mathbb{T}^2) - 2c_2(\mathbb{T}^4) = 0
\end{array}
\end{equation}

Notice that the above computation also follows easily from Theorem~\ref{thm:main1}. Since $c_{1}(X_{\psi}) = 0$ and  $c_3(X_{\psi}) = 0$, ${c_1}^3(X_{\psi})$ and ${c_1c_2}(X_{\psi})$ obviously vanishes. 

\section{Proof of Theorem~\ref{thm:main2}}
{\label{sec: FCY}}

\begin{proof} To construct a family of Calabi-Yau symplectic\/ $6$-manifolds with the fundamental group $\Z$, we again form the symplectic fiber sum of $E(1) \times \mathbb{T}^2$ with itself as in the proof of Theorem~\ref{thm:main1}. Now we choose the gluing diffemorphism $\psi$ differently. 

Let us choose the gluing diffeomorphism $\psi' : \partial(F \times \mathbb{T}^2 \times D^2) \longrightarrow  \partial(F' \times \mathbb{T}^2 \times D^2)$ that comes from an orientation preserving diffemorphism of $\mathbb{T}^{4}$ which sends the elements of $\pi_1$ as follows:
\begin{equation*}
a = 1 \mapsto b'= 1, \ \
b = 1 \mapsto d', \ \
c \mapsto c', \ \
d \mapsto a' = 1 , \ \
\mu = 1 \mapsto {\mu'}^{-1} = 1, \\
\end{equation*}

By Seifert-Van Kampen Theorem, we get the following presentation for the fundamental group of $X_{\psi'}$.

\begin{eqnarray}
\pi_1(X_{\psi'})  &=& \langle
c, d;\ c', d' \mid c=c',\, d=1,\, d'=1\rangle
\end{eqnarray}

We compute $\pi_1(X_{\psi'}) = \Z$. In particular, $X_{\psi'}$ is non-K\"{a}hler. To make $X_{\psi'}$ symplectic, we first need to perturb the ambient symplectic form in one copy to make the lagrangian tori $b \times c$ and $a \times d$ the symplectic submanifolds of $\mathbb{T}^{4}$. We refer to the Lemma 1.6 ~\cite{gompf} for the existence of such perturbation.    

Once again by the argument similar to in the proof of Theorem~\ref{thm:main1}, $c_1(X_{\psi'})$ evaluates to zero on $2$-dimensional homology classes of $X_{\psi'}$. Notice that in  this case $X_{\psi'}$ has basis for $2$-dimensional homology which in addition to the classes in the proof of Theorem~\ref{thm:main1} contains one pair of essential rim torus and a dual $2$ sphere, but again by Lemma 2.4 in \cite{IP}, $c_1(X_\psi)$ evaluates on rim tori and sphere to zero. We conclude that $X_{\psi'}$ is Calabi-Yau. 
\end{proof}

\vspace{3pt}

\begin{rmk}{\label{R1}} 
Notice that we have $SL(2, \Z) \times SL(2, \Z)$ worth of choices for our gluing map $\psi$. It was pointed to us by Robert Gompf the group action doesn't affect the resulting manifold $X_{\psi}$. By using more general gluing diffemorphism from $SL(4, \Z)$, we can also construct examples with fundamental groups $\Z_p \times \Z_q$ for $p, q \geq 1$, $\Z \times \Z_q$ or $\Z \times \Z_p$. All the symplectic CY $6$-manifolds constructed in this paper are symplectically minimal. This could be proved using \cite{LR}. The above construction can obviously be generalized to dimension $4n + 2$ for any $n \geq 2$ by using the symplectic manifolds of the form $E(1) \times \cdots \times E(1) \times \mathbb{T}^2$ and summing them along $\mathbb{T}^{2n+2}$. We will present an alternative proof of Theorem~\ref{thm:main1} and ~\ref{thm:main2} in the sequel \cite{AL}  \end{rmk}.

\begin{rmk}{\label{R1}} Using a genus $g$ Lefschetz fibration on $\CP\# (4g+5)\CPb$ over $S^2$ with the global monodromy $(a_1a_2 \cdots {a_{2g+1}}^2 \cdots a_2a_1)^2 = 1$ and forming the symplectic connected sums of $\CP\# (4g+5)\CPb \times \Sigma_g$ along $\Sigma_g \times \Sigma_g$, using a special gluing diffeomorphism that interchanges two copies of $\Sigma_g$, one can generalize Theorems \ref{thm:main1} and \ref{thm:main2} to get simply connected symplectic $6$-manifolds with Chern numbers ${c_1}^3 = 24(g-1)^2$, \ $c_1c_2 = 24(1-g)$, \ $c_3 = 8(g+2)(1-g)$. More familes can be constructed using $X(n,g) \times \Sigma_g$, $Y(n,g) \times  \Sigma_g$, and $Z(n,g) \times \Sigma_g$,
where $X(n,g)$, $Y(n,g)$ and $Z(n,g)$ are the total spaces of the $n$ fold fiber sum of three well known hyperelliptic Lefschetz fibrations given by the monodromies  $(a_1a_2 \cdots {a_{2g+1}}^2 \cdots a_2a_1)^2 = 1$, $(a_1a_2 \cdots a_{2g+1})^{2g+2} = 1$, and $(a_1a_2 \cdots a_{2g})^{4g+2} = 1$ in the mapping class group $M_g$. By considering the symplectic building blocks and ideas from \cite{gompf, Halic, A3, Ak, AP, AP2}, one can construct simply connected \emph{potentially minimal} symplectic $6$-manifolds with Chern numbers ${c_1}^3 = a$, \ $c_1c_2 = b$, \ $c_3 = c$ for any triples $(a, b, c)$ with $a\equiv 0 \bmod{2}$, $b\equiv 0 \bmod{24}$, and $c\equiv 0 \bmod{2}$. For more details and the additional building blocks needed, we refer the reader to \cite{gompf, Halic, A3, Ak, AP, ABBKP, AP2}. Unfortunately, currently there is no way to prove minimality of these examples in a general setting, but in Calabi-Yau case, one can show minimality using \cite{LR}. Also, using the building blocks from \cite{gompf, A3, Ak, AP, ABBKP, AP2}, it is simple to construct symplectic $6$-manifolds with arbitrary finitely presented group as the fundamental group and varying Chern numers $({c_1}^3, \ c_1c_2, \ c_3)$ as above. For example, if $G$ is any finitely presented group, then there exists a spin symplectic $4$-manifold $S_G$ with $\pi_1(S_G) = G$, ${c_1}^2(S_G) = 0$, and $\chi_{h}(S_G) > 0$ \cite{gompf}. $S_G$ contains a symplectic torus $\mathbb{T}^2$ of self-intersection $0$ such that the inclusion induced homomorphism $\pi_1(\mathbb{T}^2) \rightarrow \pi_1(S_G)$ is trivial. We can form the twisted sums of $\CP\# (4g+5)\CPb \times \mathbb{T}^2$ (or $X(n,g) \times \mathbb{T}^2$, $Y(n,g) \times  \mathbb{T}^2$, and $Z(n,g) \times \mathbb{T}^2$), and $S_G \times \Sigma_g$ along  $\Sigma_g \times \mathbb{T}^2$ to get symplectic $6$-manifold with the fundamental group $G$ and varying Chern numbers $({c_1}^3, \ c_1c_2, \ c_3)$. 
\end{rmk} 

\vspace{3pt}

\section{Simply Connected Symplectic CY $6$-Manifolds via Coisotropic Luttinger Surgery}{\label{CLS}} 

We would like to remark that we can also construct the symplectic Calabi-Yau $6$-manifolds as above via coisotropic Luttinger surgery. The prospect of obtaining such examples was mentioned in \cite{BK}, but no examples with $b_1 < 2$ was given. In particular, the authors were not able to obtain the simply connected examples in \cite{BK} (see page 2). We now mention our construction. We refer the reader to \cite{BK, Ho} for the definition of coisotropic Luttinger surgery. For background on Luttinger surgery, we refer to \cite{luttinger, ADK}, see also \cite{FPS, ABP, AP2, ABBKP} for the applications. 

First, we take two copies of $W = E(1) \times T^2$ and form their symplectic connected sum along the symplectic submanifolds $F \times \mathbb{T}^2$ and $F' \times \mathbb{T}^2$. We choose our gluing diffeomorphism $\phi : \partial(F \times \mathbb{T}^2 \times D^2) \longrightarrow  \partial(F' \times \mathbb{T}^2 \times D^2)$ that comes from an orientation preserving diffemorphism of $\mathbb{T}^4$ which sends the generators of $\pi_1$ as follows:
\begin{equation*}
a = 1 \mapsto a' = 1, \ \
b = 1 \mapsto b' = 1, \ \
c \mapsto c' , \ \
d \mapsto d'  , \ \
\mu = 1 \mapsto {\mu'}^{-1} = 1, \\
\end{equation*}

The resulting manifold is clearly Calabi-Yau $6$-manifold $K3 \times \mathbb{T}^2$. Next, we identify the following two $4$-tori in $K3 \times \mathbb{T}^2$: $T_1 := (a \times c) \times (d \times s)$ and $T_2 := (b \times d) \times (c \times s')$, where $s$ and $s'$ are "rim" circles of $K3$ surface, the meridians of $F^2 \times \mathbb{T}^2$. Let show that each of these $4$-tori $T_i$ has a dual $2$ sphere $S_i$ of self-intersection $-2$. Notice that each of the above mentioned $4$-dimensional torus in $\mathbb{T}^5 = \partial(E(1) \times \mathbb{T}^2 \setminus F \times \mathbb{T}^2 \times D^2)$ has a dual circle in $\mathbb{T}^5$ ("the remaining fifth circle") intersecting $T_i$ at a point, and since these circles $a$ and $b$ are null-homotopic in $\pi_1(E(1) \times \mathbb{T}^2 \setminus F \times \mathbb{T}^2 \times D^2) = \Z \times \Z$, they can be contracted in $E(1) \times \mathbb{T}^2 \setminus F \times \mathbb{T}^2 \times D^2$. The dual spheres $S_1$ and $S_2$ are obtained by contracting the circles $a$ and $b$ on both sides, using the vanishing disks of $a$ and $b$. Notice that the meridian of $T_i$ lies on $S_i$, thus null-homotopic in the in the fundamental group of complement of $\pi_1(K3 \times \mathbb{T}^2 \setminus (\nu(T_1) \cup \nu(T_2))$. Now we perform the following two coisotropic Luttinger surgeries on $4$-tori $T_1$ and $T_2$ in $K3 \times \mathbb{T}^2$: $(T_1, c^p, \pm 1)$ and $(T_2, d^q, \pm 1)$. We denote the resulting symplectic manifold by $M_{p,q}$, where $p,q \geq 0$. Using the dual $-2$ spheres $S_i$ of $T_{i}$ and the meridian of $T_i$ are null-homotopic in the in the fundamental group $\pi_1(K3 \times \mathbb{T}^2 \setminus (\nu(T_1) \cup \nu(T_2))$,  we easily see that the fundamental groups of $M_{p,q}$ are one of the folowing abelian groups:  $\Z_p \times \Z_q$ for $p, q \geq 1$, $\Z \times \Z_q$ or $\Z \times \Z_p$ for $p = 0$, $q \geq 1$ or  $p \geq 1$, $q = 0$. The first coisotropic Luttinger surgery gives $c^p = 1$ and the second surgery produces the relation $d^q = 1$ in $\pi_1(K3 \times \mathbb{T}^2) = \Z \times \Z$. If we set $p = q = 1$, then $M_{1,1}$ has a trivial fundamental group. By setting $p = 1$, $q = 0$ or $p = 0$, $q = 1$, we get the symplectic $6$-manifolds $M_{1,0}$ and $M_{0,1}$ with fundamental group $\Z$. 

Let us now prove that the simply connected symplectic manifolds $M_{1,1}$ is Calabi-Yau. Other cases are similar, and we study them in \cite{AL}.

To get a basis for $H_2$ of $K3 \times \mathbb{T}^2$ ($b_2 = 23$), we use $22$ two dimensional classes in Example ~\ref{E} plus $pt \times \mathbb{T}^2$. When we do Luttinger surgery to get the symplectic $6$-manifolds $M_{1,1}$, we kill $2$ rim 4-tori and their dual $2$ spheres in homology, thus $c_1(M_{1,1})$ automatically zero on rim $2$-tori and these sphere classes. $c_1(M_{1,1})$ is zero on the classes coming from two copies of Milnor fiber $\Phi(1) \in K3 \times \mathbb{T}^2$. By Lemma 2.4 in \cite{IP}, $c_1(M_{1,1})$ is also zero on sphere section $\sigma$ as it is the union of $-1$ sphere sections of $E(1)$ and the remaining classes $f$ and $pt \times \mathbb{T}^2$. To apply Lemma 2.4 in \cite{IP} in this setting, one first has to view $K3 \times \mathbb{T}^2$ as the fiber sums of $E(1) \times \mathbb{T}^2$. We also refer the reader to \cite{HoLi}, where coisotropic Luttinger surgery was suggested as candidate for a symplectic Calabi-Yau surgery.

\begin{rmk} Similar to our construction as above, we can obtain simply connected symplectic $6$-manifolds by performing two coisotropic Luttinger surgeries on two $4$-tori in $W(n) = E(n) \times \mathbb{T}^2$ (for $n \geq 3$) and  $W(n, K) = E(n)_{K} \times \mathbb{T}^2$ (for $n \geq 2$), where $E(n)_{K}$ is a symplectic homotopy elliptic surface of Fintushel and Stern \cite{FS2}. We would like to remark that the choice of these $4$-tori are same as in the proof of Theorems~\ref{thm:main1}. These examples are not Calabi-Yau. One can also obtain the simply connected symplectic $6$-manifolds by performing $2g$ coisotropic Luttinger surgeries on $X(n, g) \times \Sigma_g$, $Y(n, g) \times \Sigma_g$, $Z(n, g) \times \Sigma_g$, where $X(n,g)$, $Y(n,g)$ and $Z(n,g)$ are the total spaces of the $n$ fold fiber sum of three well known hyperelliptic Lefschetz fibrations given by the monodromies  $(a_1a_2 \cdots {a_{2g+1}}^2 \cdots a_2a_1)^2 = 1$, $(a_1a_2 \cdots a_{2g+1})^{2g+2} = 1$, and $(a_1a_2 \cdots a_{2g})^{4g+2} = 1$ in the mapping class group $M_g$. The abelian and cyclic $\pi_1$ examples can be obtained from $V_{m,n,g} \times \Sigma_g$, where $V_{m,n,g}$ is $m$ fold fiber sum of the genus $2g+n-1$ Lefschetz fibration on $\Sigma_g \times \mathbb{S}^2 \# 4n\CPb$ \cite{Korkmaz, Mat}, and $Sym^2(\Sigma_g) \times \mathbb{T}^2$. Notice that $V_{1,1,1}$ is well known Matsumoto's genus two fibration on $\mathbb{T}^2 \times \mathbb{S}^2 \# 4\CPb$. The details of these construction given in \cite{AL}.  

\end{rmk}

\begin{rmk}

We expect that the manifolds $M_{1,1}$ and $M_{1,0}$ are closely related (perhaps the same) to $X_{\psi}$ and $X_{\psi'}$ that we constructed in Theorems~\ref{thm:main1} and ~\ref{thm:main2}. In a forthcoming paper, we study these manifolds in details. We also use the symplectic Calabi-Yau $6$-manifolds $X_{\psi}$ and $X_{\psi'}$ along with the symplectic building blocks that we constructed above and $Sym^3(\Sigma_g)$, $3$-fold symmetric product of the genus $g$ surface, to study the geography of symplectic $6$-manifolds in \cite{AL}.

\end{rmk}

\section*{Acknowledgments} The author is grateful to Tian-Jun Li for valuable discussions and for his kind encouragement. The author is thankful to Robert Gompf for making various suggestions which have helped him in improving this paper. The author is also grateful for the support by NSF grants FRG-1065955 and DMS-1005741.


\begin{thebibliography}{99}

\bibitem{A3}  Akhmedov, A, \textit{Construction of symplectic cohomology $S^{2}\times S^{2}$}, G\"{o}kova Geometry and Topology Proceedings, {\bf 14} (2007), 36--48.

\bibitem{Ak} A. Akhmedov, 
\textit{Small exotic\/ $4$-manifolds}, 
Algebr. Geom. Topol. \textbf{8} (2008), 1781--1794. 

\bibitem{AP}  A. Akhmedov and B. D. Park,
\textit{Exotic smooth structures on small\/ $4$-manifolds}, 
Invent. Math. \textbf{173} (2008), 209--223.

\bibitem{AL}  A. Akhmedov, T. J. Li and B. D. Park,
\textit{The geography of symplectic CY 6-manifolds}, preprint.

\bibitem{ABP} A. Akhmedov, R. \.{I}. Baykur and B. D. Park, 
\textit{Constructing infinitely many smooth structures on small\/ $4$-manifolds}, 
J. Topol. \textbf{1} (2008), 409--428.

\bibitem{ABBKP}  Akhmedov, A., Baldridge, S., Baykur, I., Kirk, P., Park, B.~D, \textit{Simply connected minimal symplectic 4-manifolds with signature less than -1}, Journal of European Math Society, {\bf 1}  (2010), 133--161. 

\bibitem{AP2}  A. Akhmedov, and B. D. Park, \textit{Exotic smooth structures on small $4$-manifolds with odd signatures}, Inventiones Mathematicae, {\bf 181} (2010), 209--223.

\bibitem{ADK}  D. Auroux, S. K. Donaldson and L. Katzarkov,
\textit{Luttinger surgery along Lagrangian tori and non-isotopy
for singular symplectic plane curves}, Math. Ann. \textbf{326} (2003), 185--203.


\bibitem{B} S. Bauer, \textit{Almost complex 4-manifolds with vanishing first Chern class}, J. Diff. Geom. \textbf{79} (2008), 25--32.

\bibitem{BK} S. Baldridge and P. Kirk, \textit{Coisotropic Luttinger surgery and some new symplectic 6-manifolds}, preprint (arXiv:1105.3519v1), (2011).

\bibitem{FP1}  J. Fine and D. Panov,
\textit{Symplectic Calabi-Yau manifolds, minimal surfaces and the hyperbolic geometry of the conifold}, J. Diff. Geom. \textbf{82} (2009), 155--205.

\bibitem{FP2}  J. Fine and D. Panov,
\textit{Hyperbolic geometry and non-Kahler manifolds with trivial canonical bundle},
Geom. Topol \textbf{14} (2010), 1731--1763.

\bibitem{FP2}  J. Fine and D. Panov,
\textit{Building symplectic manifolds using hyperbolic geometry},
Proc. Gokova. Geom. Topol \textbf{17} (2009), 124--136.

\bibitem{FS2}  Fintushel, R., Stern, R. J.: 
Knots, links and $4$-manifolds.  Invent. Math. 
{\bf 134}, (1998), 363--400. 

\bibitem{FPS}  R. Fintushel, B. D. Park and R. J. Stern,
\textit{Reverse engineering small\/ $4$-manifolds}, 
Algebr. Geom. Topol. \textbf{7} (2007), 2103--2116.

\bibitem{Halic} M. Halic, \textit{On the geography of symplectic 6-manifolds}, Manuscripta Math. \textbf{99} (1999), 371--381.

\bibitem{Ho}  C. I. Ho,
\textit{Topological Methods in Symplectic Geometry}, Ph.D Thesis (2011).

 \bibitem{HoLi} C. I. Ho and T. J. Li, \textit{Luttinger surgery and Kodaira dimension}, to appear.

\bibitem{gompf} R. E. Gompf,
\textit{A new construction of symplectic manifolds},
Ann. of Math. \textbf{142} (1995), 527--595.

\bibitem{GS}  R. E. Gompf and A. I. Stipsicz, 
\textit{$4$-Manifolds and Kirby Calculus}, Graduate Studies in
Mathematics, vol. 20, Amer. Math. Soc., Providence, RI, 1999.

\bibitem{Korkmaz} M. Korkmaz, \textit{Noncomplex smooth $4$-manifolds with Lefschetz fibrations}, Internat. Math. Res. Not. 2001, No. 3, 115--128. 

\bibitem{luttinger}  K. M. Luttinger,  \textit{Lagrangian tori in\/ $\R^4$},
J. Differential Geom. \textbf{42} (1995), 220--228.

\bibitem{Mat} Y. Matsumoto, \textit{Lefschetz fibrations of genus two - a topological approach}, Proceedings of the 37th Taniguchi Symposium on “Topology and Teichm¨uller Spaces”, World Scientific,
Singapore, 1996, 123--148.

\bibitem{MS} J. Morgan and Z.  Szabo, \textit{Homotopy K3 surfaces and mod 2 Seiberg-Witten invariant}, Math. Res. Lett. \textbf{4}, no. 1, (1997), 17--21.

\bibitem{STY}  I. Smith, R. Thomas, and S.-T. Yau,
\textit{Symplectic Conifold Transitions},
J. Diff. Geom. \textbf{62} (2002), 209--242.

\bibitem{IP4}  E-N. Ionel and  T. Parker,  
\textit{Relative Gromov-Witten invariants},
 Ann. of Math. (2)  \textbf{157}, (2003), no. 1, 45--96.
 
 \bibitem{IP}  E-N. Ionel and  T. Parker,  
\textit{The symplectic sum formula for Gromov-Witten invariants},
 Ann. of Math. (2)  \textbf{159}, (2004), no. 3, 935--1025.

 \bibitem{Li1} T. J. Li, \textit{Quaternionic bundles and Betti numbers of symplectic 4-manifolds with Kodaira dimension zero}, Int. Math. Res. Not. 2006, Art. ID 37385, 28 pp.

 \bibitem{Li2} T. J. Li, \textit{Symplectic Calabi-Yau Surfaces}, Adv. Lec. in Math., \textbf{14}, Handbook of Geometric Analysis, no. 3, 2010. 

 \bibitem{Li3} T. J. Li, \textit{Symplectic 4-manifolds with Kodaira dimension zero}, J. Diff. Geom. \textbf{74}, (2006), no. (2), 321--352.

 \bibitem{LR} T. J. Li and Y. Ruan, \textit{Symplectic Birational Geometry}, preprint, arXiv:0906.3265v1, (2009).

\bibitem{Usher} M. Usher, \textit{Kodaira dimension and symplectic sums}, Commentarii Mathematici Helvetici, \textbf{84} (2009), no. 1, 57--85.

\end{thebibliography}
\end{document}